# Some power series extensions of Fibonacci and Lucas polynomials

Johann Cigler

**Abstract**

We study formal power series which can be interpreted as interpolations of Fibonacci and Lucas polynomials with even (or odd) indices.

## 1. Introduction

We consider the Fibonacci polynomials $F_n(x)$ defined by

(1) $$F_n(x) = xF_{n-1}(x) + F_{n-2}(x)$$

with initial values $F_0(x) = 0$ and $F_1(x) = 1$ and the corresponding Lucas polynomials $L_n(x) = F_{n+1}(x) + F_{n-1}(x)$ which satisfy

(2) $$L_n(x) = xL_{n-1}(x) + L_{n-2}(x)$$

with initial values $L_0(x) = 2$ and $L_1(x) = x$.

The first terms are

$$\left(F_n(x)\right)_{n \geq 0} = \left(0, 1, x, 1+x^2, 2x+x^3, 1+3x^2+x^4, 3x+4x^3+x^5, 1+6x^2+5x^4+x^6, \cdots\right),$$

$$\left(L_n(x)\right)_{n \geq 0} = \left(2, x, 2+x^2, 3x+x^3, 2+4x^2+x^4, 5x+5x^3+x^5, 2+9x^2+6x^4+x^6, \cdots\right).$$

The well-known formulae

$$F_n(x) = \sum_{k=0}^{\left\lfloor \frac{n-1}{2} \right\rfloor} \binom{n-1-k}{k} x^{n-1-2k} \text{ and } L_n(x) = \sum_{k=0}^{\left\lfloor \frac{n}{2} \right\rfloor} \binom{n-k}{k} \frac{n}{n-k} x^{n-2k}$$

imply

$$F_{2n}(x) = \sum_{k=1}^{n} \binom{n+k-1}{2k-1} x^{2k-1}, \quad F_{2n+1}(x) = \sum_{k=0}^{n} \binom{n+k}{2k} x^{2k},$$

$$L_{2n}(x) = \sum_{k=0}^{n} \binom{n+k}{2k} \frac{2n}{n+k} x^{2k}, \quad L_{2n+1}(x) = \sum_{k=0}^{n} \binom{n+k+1}{2k+1} \frac{2n+1}{n+k+1} x^{2k+1}.$$

These formulae suggest the following formal power series depending on a parameter $t \in \mathbb{R}$ which can be interpreted as interpolations of these polynomials.

(3)
$$\Phi_0(2t, x) = \sum_{k=1}^{\infty} \binom{t+k-1}{2k-1} x^{2k-1},$$
$$\Phi_1(2t+1, x) = \sum_{k=0}^{\infty} \binom{t+k}{2k} x^{2k},$$
$$\Lambda_0(2t, x) = \sum_{k=0}^{\infty} \binom{t+k}{2k} \frac{2t}{t+k} x^{2k},$$
$$\Lambda_1(2t+1, x) = \sum_{k=0}^{\infty} \binom{t+k+1}{2k+1} \frac{2t+1}{t+k+1} x^{2k+1}.$$

## 2. Closed formulae for these series

We are looking for closed formulas for these series which generalize Binet's formulas

(4)
$$F_n(x) = \frac{\alpha^n - \bar{\alpha}^n}{\alpha - \bar{\alpha}}$$
$$L_n(x) = \alpha^n + \bar{\alpha}^n$$

where

(5)
$$\alpha = \alpha(x) = \frac{x + \sqrt{x^2+4}}{2},$$
$$\bar{\alpha} = \bar{\alpha}(x) = \frac{-1}{\alpha(x)} = \frac{x - \sqrt{x^2+4}}{2}$$

are the roots of the characteristic polynomial $z^2 - xz - 1$ of the recurrence of the Fibonacci and Lucas polynomials.

Note that

(6)
$$F_{2n}(x) = \frac{\alpha^{2n} - \alpha^{-2n}}{\sqrt{x^2+4}},$$
$$F_{2n+1}(x) = \frac{\alpha^{2n+1} + \alpha^{-2n-1}}{\sqrt{x^2+4}},$$
$$L_{2n}(x) = \alpha^{2n} + \alpha^{-2n},$$
$$L_{2n+1}(x) = \alpha^{2n+1} - \alpha^{-2n-1}.$$

The identity

$$\alpha(e^y - e^{-y}) = \frac{e^y - e^{-y} + \sqrt{(e^y - e^{-y})^2 + 4}}{2} = \frac{e^y - e^{-y} + e^y + e^{-y}}{2} = e^y$$

can be written as

(7) $$\alpha(2\sinh(y)) = e^y$$

where $\sinh(y) = \dfrac{e^y - e^{-y}}{2}$ denotes the hyperbolic sine series.

Let $\sinh^{-1}(x)$ be the formal reverse of $\sinh(x)$ which satisfies $\sinh^{-1}(\sinh(x)) = \sinh(\sinh^{-1}(x)) = x$.

By setting $y = \sinh^{-1}\left(\dfrac{x}{2}\right)$ we get from (7)

(8) $$\alpha(x) = e^{\sinh^{-1}\left(\frac{x}{2}\right)}.$$

It is easily verified by induction that

(9) $$\left(x + \dfrac{1}{x}\right) F_{2n}\left(x - \dfrac{1}{x}\right) = x^{2n} - \dfrac{1}{x^{2n}}$$

holds. Therefore, we get for $x = e^y$

$(e^y + e^{-y}) F_{2n}(e^y - e^{-y}) = e^{2ny} - e^{-2ny}$ which gives

(10) $$(e^y + e^{-y}) \sum_{k=1}^{n} \binom{n+k-1}{2k-1} (e^y - e^{-y})^{2k-1} = e^{2ny} - e^{-2ny}.$$

This implies that for all $t \in \mathbb{R}$

(11) $$(e^y + e^{-y}) \sum_{k=1}^{\infty} \binom{t+k-1}{2k-1} (e^y - e^{-y})^{2k-1} = e^{2ty} - e^{-2ty},$$

because for $k \in \mathbb{N}$ the coefficient of $y^k$ of each side is a polynomial in $t$. Since by (10) the two polynomials are equal for infinitely many $t = n \in \mathbb{N}$ they are equal as polynomials.

Identity (11) can be written as

(12) $$\sqrt{(e^y - e^{-y})^2 + 4} \, \Phi_0(2t, e^y - e^{-y}) = (e^y + e^{-y}) \Phi_0(2t, e^y - e^{-y}) = e^{2ty} - e^{-2ty}.$$

By setting $x = e^y - e^{-y} = 2\sinh(y)$ we get $\sqrt{x^2 + 4} \, \Phi_0(2t, x) = 2\sinh\left(2t \sinh^{-1}\left(\dfrac{x}{2}\right)\right)$.

Thus, we get by (8)

(13) $$\Phi_0(t, x) = \dfrac{2 \sinh\left(t \sinh^{-1}\left(\dfrac{x}{2}\right)\right)}{\sqrt{x^2 + 4}} = \dfrac{\alpha(x)^t - \alpha(x)^{-t}}{\sqrt{x^2 + 4}}.$$

In the same way

$$(14) \quad \left(x+\frac{1}{x}\right)F_{2n+1}\left(x-\frac{1}{x}\right) = x^{2n+1} + \frac{1}{x^{2n+1}}$$

gives

$$(15) \quad \Phi_1(t,x) = \frac{2\cosh\left(t\sinh^{-1}\left(\frac{x}{2}\right)\right)}{\sqrt{x^2+4}} = \frac{\alpha(x)^t + \alpha(x)^{-t}}{\sqrt{x^2+4}}.$$

Similarly

$$(16) \quad L_{2n}\left(x-\frac{1}{x}\right) = x^{2n} + \frac{1}{x^{2n}}$$

implies

$$\Lambda_0(t,x) = 2\cosh\left(t\sinh^{-1}\left(\frac{x}{2}\right)\right) = \alpha(x)^t + \alpha(x)^{-t}$$

and

$$(17) \quad L_{2n+1}\left(x-\frac{1}{x}\right) = x^{2n+1} - \frac{1}{x^{2n+1}}$$

gives

$$(18) \quad \Lambda_1(t,x) = 2\sinh\left(t\sinh^{-1}\left(\frac{x}{2}\right)\right) = \alpha(x)^t - \alpha(x)^{-t}.$$

In order to state the results efficiently, we interpret the index $j$ of $\Phi_j(t,x)$ and $\Lambda_j(t,x)$ modulo 2 and can write

$$(19) \quad \Phi_j(t,x) = \frac{\alpha(x)^t - (-1)^j \alpha(x)^{-t}}{\sqrt{x^2+4}}, \quad \Lambda_j(t,x) = \alpha(x)^t + (-1)^j \alpha(x)^{-t}.$$

This gives

$$(20) \quad \frac{\Lambda_j(t,x) + \sqrt{x^2+4}\,\Phi_j(t,x)}{2} = \alpha(x)^t$$

and comparing with (3)

$$(21) \quad \alpha(x)^t = \frac{\Lambda_0(t,x) + \Lambda_1(t,x)}{2} = \sum_{n \geq 0} \binom{\frac{n+t}{2}}{n} \frac{t}{n+t} x^n.$$

**Remark**

It should be noted that Sergio Falcon and Angel Plaza [1] studied the same situation from another point of view. As analogs of Binet's formula for $F_n(k)$ they considered the continuous functions

$$\frac{\alpha(k)^t - (-1)^j \alpha(k)^{-t}}{\sqrt{k^2 + 4}}$$ and studied their algebraic and geometric properties.

## 3. Some simple identities

Generalizing $L_n(x) = F_{n-1}(x) + F_{n+1}(x)$ we get using (19) and $\alpha(x) + \alpha(x)^{-1} = \alpha - \bar{\alpha} = \sqrt{x^2 + 4}$

(22) $$\Lambda_j(t, x) = \Phi_{j+1}(t-1, x) + \Phi_{j+1}(t+1, x).$$

From $\alpha^2 = x\alpha + 1$ and $\alpha^{-2} = -x\alpha^{-1} + 1$ we see that

$$\alpha(x)^{t+2} \pm \alpha(x)^{-t-2} = x\left(\alpha(x)^{t+1} \mp \alpha(x)^{-t-1}\right) + \alpha(x)^t \pm \alpha(x)^{-t}.$$

This gives

(23) $$\begin{aligned}\Phi_j(t+2, x) &= x\Phi_{j+1}(t+1, x) + \Phi_j(t, x), \\ \Lambda_j(t+2, x) &= x\Lambda_{j+1}(t+1, x) + \Lambda_j(t, x)\end{aligned}$$

generalizing the recurrence of the Fibonacci and Lucas polynomials.

Since

$$\left(\alpha(x)^{t+1} \pm \alpha(x)^{-t-1}\right)\left(\alpha(x)^{t-1} \pm \alpha(x)^{-t+1}\right) - \left(\alpha(x)^t \mp \alpha(x)^{-t}\right)^2$$
$$= \alpha(x)^{2t} \pm \alpha(x)^{-2} \pm \alpha(x)^2 + \alpha(x)^{-2t} - \alpha(x)^{2t} - \alpha(x)^{-2t} \pm 2 = \pm\left(\alpha(x)^2 + \alpha(x)^{-2} + 2\right) = \pm\left(x^2 + 4\right)$$

we get

(24) $$\begin{aligned}\Phi_j(t+1, x)\Phi_j(t-1, x) - \Phi_{j+1}(t, x)^2 &= (-1)^{j+1}, \\ \Lambda_j(t+1, x)\Lambda_j(t-1, x) - \Lambda_{j+1}(t, x)^2 &= (-1)^j\left(x^2 + 4\right).\end{aligned}$$

For $t = 2n$ and $j = 1$ or $t = 2n+1$ and $j = 0$ these reduce to Cassini's formula $F_{n+1}(x)F_{n-1}(x) - F_n^2(x) = (-1)^n$ and to $L_{n+1}(x)L_{n-1}(x) - L_n^2(x) = (-1)^{n-1}\left(x^2 + 4\right).$

From

$$\left(\alpha(x)^{t+1} \pm \alpha(x)^{-t-1}\right)\left(\alpha(x)^{t-1} \pm \alpha(x)^{-t+1}\right) - \left(\alpha(x)^t \pm \alpha(x)^{-t}\right)^2$$
$$= \alpha(x)^{2t} \pm \alpha(x)^{-2} \pm \alpha(x)^2 + \alpha(x)^{-2t} - \alpha(x)^{2t} - \alpha(x)^{-2t} \mp 2 = \pm\left(\alpha(x)^2 + \alpha(x)^{-2} - 2\right) = \pm x^2$$

we get

(25)
$$\Phi_j(t+1,x)\Phi_j(t-1,x) - \Phi_j(t,x)^2 = (-1)^{j+1}\frac{x^2}{x^2+4},$$
$$\Lambda_j(t+1,x)\Lambda_j(t-1,x) - \Lambda_j(t,x)^2 = (-1)^j x^2.$$

Let us also note that (19) implies

(26)
$$\Phi_0(2t,x) = \Phi_0(t,x)\Lambda_0(t,x),$$
$$\Lambda_0(2t,x) = \Lambda_0^2(t,x) - 2,$$

which generalize $F_{2n}(x) = F_n(x)L_n(x)$ and $L_{2n}(x) + 2 = L_n^2(x)$.

For positive integers $n$ we have $\Phi_j(2n+j,x) = F_{2n+j}(x)$, $\Lambda_j(2n+j,x) = L_{2n+j}(x)$ by definition and for the other integer values we get

(27)
$$\Lambda_{j+1}(2n+j,x) = F_{2n+j}(x)\sqrt{x^2+4},$$
$$\Phi_{j+1}(2n+j,x) = \frac{L_{2n+j}(x)}{\sqrt{x^2+4}},$$

because the first line follows from

$$\Lambda_{j+1}(2n+j,x) = \alpha(x)^{2n+j} + \frac{(-1)^{j+1}}{\alpha(x)^{2n+j}} = \sqrt{x^2+4}\frac{1}{\sqrt{x^2+4}}\left(\alpha(x)^{2n+j} + \frac{(-1)^{j-1}}{\alpha(x)^{2n+j}}\right)$$
$$= \sqrt{x^2+4}\Phi_j(2n+j,x) = \sqrt{x^2+4}F_{2n+j}(x)$$

and from (20) we then get
$$\Lambda_{j+1}(2n+j,x) + \sqrt{x^2+4}\Phi_{j+1}(2n+j,x) = 2\alpha^{2n+j} = L_{2n+j}(x) + \sqrt{x^2+4}F_{2n+j}(x)$$ which gives

$$\Phi_{j+1}(2n+j,x) = \frac{L_{2n+j}(x)}{\sqrt{x^2+4}}.$$

Finally let us consider in detail the values $\Phi_j(k) := \Phi_j(k,1)$ and $\Lambda_j(k) := \Lambda_j(k,1)$ for positive integers $k$ and their relation to the Fibonacci numbers $F_k$ and the Lucas numbers $L_k$.

| $F_k$ | 0 | 1 | 1 | 2 | 3 | 5 | 8 | 13 | 21 |
|---|---|---|---|---|---|---|---|---|---|
| $\Phi_0(k)$ | 0 | $\frac{1}{\sqrt{5}}$ | 1 | $\frac{4}{\sqrt{5}}$ | 3 | $\frac{11}{\sqrt{5}}$ | 8 | $\frac{29}{\sqrt{5}}$ | 21 |
| $\Phi_1(k)$ | $\frac{2}{\sqrt{5}}$ | 1 | $\frac{3}{\sqrt{5}}$ | 2 | $\frac{7}{\sqrt{5}}$ | 5 | $\frac{18}{\sqrt{5}}$ | 13 | $\frac{47}{\sqrt{5}}$ |

They satisfy recurrence (23) and can also be obtained from the Fibonacci numbers by using (27) which gives $\Phi_{j+1}(2n+j) = \dfrac{F_{2n+j+1} + F_{2n+j+1}}{\sqrt{5}}$ or by using (25) which gives

$$\Phi_0(2k+1) = \frac{\sqrt{1+5F_{2k}F_{2k+2}}}{\sqrt{5}} \text{ and } \Phi_1(2k) = \frac{\sqrt{5F_{2k+1}F_{2k-1}-1}}{\sqrt{5}}.$$

So, for example we get $\Phi_0(5) = \frac{3+8}{\sqrt{5}} = \frac{\sqrt{5\times 3\times 8+1}}{\sqrt{5}}$, $\Phi_1(4) = \frac{2+5}{\sqrt{5}} = \frac{\sqrt{5\times 2\times 5-1}}{\sqrt{5}}$.

Similarly we get for the Lucas numbers

| $L_k$ | 2 | 1 | 3 | 4 | 7 | 11 | 18 | 29 | 47 |
|---|---|---|---|---|---|---|---|---|---|
| $\Lambda_0(k)$ | 2 | $\sqrt{5}$ | 3 | $2\sqrt{5}$ | 7 | $5\sqrt{5}$ | 18 | $13\sqrt{5}$ | 47 |
| $\Lambda_1(k)$ | 0 | 1 | $\sqrt{5}$ | 4 | $3\sqrt{5}$ | 11 | $8\sqrt{5}$ | 29 | $21\sqrt{5}$ |

From $\Lambda_0(2k) = L_{2k}$ we get from (25) $\Lambda_0(2k+1) = \sqrt{L_{2k}L_{2k+2}-1}$ and $\Lambda_1(2k) = \sqrt{L_{2k-1}L_{2k+1}+1}$.

For example $\Lambda_0(3) = \frac{3+7}{\sqrt{5}} = \sqrt{3\times 7-1}$, $\Lambda_1(4) = \frac{4+11}{\sqrt{5}} = \sqrt{4\times 11+1}$.

Let us also note that $\Lambda_j(k) + \sqrt{5}\Phi_j(k) = 2\alpha(1)^k = 2\left(\frac{1+\sqrt{5}}{2}\right)^k = L_k + \sqrt{5}F_k$.

**Remark**

It should be noted that R. Witula [2] studied the complex-valued functions
$F_t = \frac{\alpha(1)^t - e^{i\pi t}\alpha(1)^{-t}}{\sqrt{5}}$ and $L_t = \alpha(1)^t + e^{i\pi t}\alpha(1)^{-t}$ of the real variable $t$ which reduce to the Fibonacci and Lucas numbers for all integers $t = n$.

**References**

[1] Sergio Falcon and Angel Plaza, The k-Fibonacci hyperbolic functions, Chaos, Solitons and Fractals 38(2008), 409-420

[2] R. Witula, Fibonacci and Lucas Numbers for Real Indices and Some Applications, Acta Physica Polonica A 120 (2011), 755-758